\documentclass[10pt]{article}
\usepackage{mathpazo}
\usepackage{amsmath}
\usepackage{amsthm}
\usepackage{amssymb}
\usepackage{url}
\usepackage{latexsym}
\usepackage[xcolor=pst]{pstricks}
\usepackage{graphicx, pst-plot, pst-node, pst-text, pst-tree}
\usepackage{titlefoot}
\usepackage[small]{titlesec}
\usepackage{units} 
\usepackage[small,it]{caption}

\usepackage{color}
\definecolor{lightgray}{rgb}{0.8, 0.8, 0.8}
\definecolor{darkgray}{rgb}{0.7, 0.7, 0.7}
\definecolor{darkblue}{rgb}{0, 0, .4}

\usepackage[bookmarks]{hyperref}
\hypersetup{
        colorlinks=true,
        linkcolor=darkblue,
        anchorcolor=darkblue,
        citecolor=darkblue,
        urlcolor=darkblue,
        pdfpagemode=UseThumbs,
        pdftitle={Generating and Enumerating $321$-Avoiding and Skew-Merged Simple Permutations},
        pdfsubject={Combinatorics},
        pdfauthor={Albert and Vatter},
}

\newcounter{todocounter}

\theoremstyle{plain}

\theoremstyle{definition}

\makeatletter
\newfont{\footsc}{cmcsc10 at 8truept}
\newfont{\footbf}{cmbx10 at 8truept}
\newfont{\footrm}{cmr10 at 10truept}
\renewenvironment{abstract}%
                {
                  \begin{list}{}%
                     {\setlength{\rightmargin}{1in}%
                      \setlength{\leftmargin}{1in}}%
                   \item[]\ignorespaces\begin{small}}%
                 {\end{small}\unskip\end{list}}

\newcommand{\C}{\mathcal{C}}

\renewcommand{\S}{\mathcal{S}}

\newcommand{\OEISlink}[1]{\href{http://oeis.org/#1}{#1}}
\newcommand{\OEISref}{\href{http://oeis.org/}{OEIS}~\cite{sloane:the-on-line-enc:}}
\newcommand{\OEIS}[1]{(Sequence \OEISlink{#1} in the \OEISref.)}

\newcommand{\Geom}{\operatorname{Geom}}

\newcommand{\st}{\::\:}

\amssubj{05A05, 05A15}

\newpagestyle{main}[\small]{
        \headrule
        \sethead[\usepage][][]
        {\sc Generating and Enumerating $321$-Avoiding and Skew-Merged Simple Permutations}{}{\usepage}}

\title{\sc Generating and Enumerating $321$-Avoiding and Skew-Merged Simple Permutations}
\author{%
Michael H. Albert\footnote{This project is the result of a research visit of the authors to the University of St Andrews, which was supported by EPSRC via the grant EP/J006440/1.}\\[-0.25ex]
\small Department of Computer Science\\[-0.5ex]
\small University of Otago\\[-0.5ex]
\small Dunedin, New Zealand\\[1.5ex]
Vincent Vatter\footnotemark[\value{footnote}]\footnote{Vatter was also partially supported by the NSA Young Investigator Grant H98230-12-1-0207.}\\[-0.25ex]
\small Department of Mathematics\\[-0.5ex]
\small University of Florida\\[-0.5ex]
\small Gainesville, Florida USA\\[-1.5ex]
}

\titleformat{\section}
        {\large\sc}
        {\thesection.}{1em}{}

\date{}

\begin{document}
\maketitle

\pagestyle{main}

\begin{abstract}
The simple permutations in two permutation classes --- the $321$-avoiding permutations and the skew-merged permutations --- are enumerated using a uniform method.  In both cases, these enumerations were known implicitly, by working backwards from the enumeration of the class, but the simple permutations had not been enumerated explicitly.  In particular, the enumeration of the simple skew-merged permutations leads to the first truly structural enumeration of this class as a whole.  The extension of this method to a wider collection of classes namely grid classes of infinite paths is discussed.
\end{abstract}

\section{Introduction}


Given permutations $\pi$ and $\sigma$, thought of as sequences of positive integers, we say that $\pi$ \emph{contains} $\sigma$, and write $\sigma\le\pi$, if $\pi$ has a subsequence $\pi(i_1)\cdots\pi(i_k)$ of the same length as $\sigma$ which is order isomorphic to $\sigma$ (i.e.~$\pi(i_s) < \pi(i_t)$ if and only if $\sigma(s) < \sigma(t)$); otherwise, we say that $\pi$ \emph{avoids} $\sigma$.  Containment is a partial order on permutations, and is the only such order considered in this paper. For example, $\pi=391867452$ contains $\sigma=51342$, as can be seen by considering the subsequence $\pi(2)\pi(3)\pi(5)\pi(6)\pi(9)=91672$.

Given a collection $X$ of permutations, we denote by $X_n$ the set of permutations in $X$ of length $n$.  The \emph{generating function} of $X$ is then
$$
\sum_{\pi\in X} x^{|\pi|}=\sum_{n}|X_n|x^n,
$$
where $|\pi|$ denotes the length of $\pi$.  As a matter of convention, except when explicitly stated otherwise, we do not include the empty permutation in our generating functions.  We are particularly interested in collections of permutations $\C$ which are downwards closed in the containment order, i.e.~$\pi\in\C$ and $\sigma\le\pi$, then $\sigma\in\C$; we call such collections \emph{permutation classes}.

The two permutation classes we consider in this paper are the $321$-avoiding permutations and the \emph{skew-merged} permutations.  This latter class is defined as the set of permutations  which can be written as the union of an increasing and a decreasing subsequence.  Stankova~\cite{stankova:forbidden-subse:} proved that the skew-merged permutations can also be characterized as the permutations that avoid both $2143$ and $3412$.  This class was enumerated by Atkinson~\cite{atkinson:permutations-wh:} via a rather intricate argument.

Our interest is with the simple permutations in these classes, and to discuss these we need a few preliminary definitions.  An \emph{interval} in the permutation $\pi$ is a set of contiguous indices $I=\{a,a+1,\dots,b\}$ such that the set $\{\pi(i)\st i\in I\}$ is also contiguous.  Every permutation $\pi$ of length $n$ has \emph{trivial intervals} of lengths $0$, $1$, and $n$, and other intervals are called \emph{proper}.  A permutation is called \emph{simple} if it does not have any proper intervals.

Simple permutations are precisely those that do not arise from a non-trivial inflation, in the following sense.  Given a permutation $\sigma$ of length $m$ and nonempty permutations $\alpha_1,\dots,\alpha_m$, the \emph{inflation} of $\sigma$ by $\alpha_1,\dots,\alpha_m$,  denoted $\sigma[\alpha_1,\dots,\alpha_m]$, is the permutation of length $|\alpha_1|+\cdots+|\alpha_m|$ obtained by replacing each entry $\sigma(i)$ by an interval that is order isomorphic to $\alpha_i$ in such a way that the intervals are order isomorphic to $\sigma$.  For example,
\[
2413[1,132,321,12]=4\ 798\ 321\ 56.
\]
It can be established that every permutation is the inflation of a unique simple permutation, called its \emph{simple quotient} and, moreover, that the intervals in such an inflation are unique unless the simple quotient is $12$ or $21$. Permutations whose simple quotient is $12$ are called \emph{sum decomposable} and those whose simple quotient is $21$ are called \emph{skew decomposable}. For enumerative purposes we specify a unique representation of sum decomposable permutations, expressing the sum decomposable $\pi$ as $\pi=12[\alpha,\beta]$ where $\alpha$ is sum indecomposable. We represent skew decomposable permutations uniquely in an analogous way as $21[\alpha, \beta]$ where $\alpha$ is skew indecomposable.  Inflations of these two permutations occur frequently enough that we give them special notation, writing $\alpha\oplus\beta$ for $12[\alpha,\beta]$ and $\alpha\ominus\beta$ for $21[\alpha,\beta]$.

\section{Implicit Derivation}
\label{sec-implicit-321}

Before introducing our method of counting simple $321$-avoiding permutations directly, we show how this enumeration can be obtained through an implicit relation between it and the enumeration of the full class.  The nonempty $321$-avoiding permutations have the generating function
\[
c(x) = \frac{1-2x-\sqrt{1-4x}}{2x}
=
x+2x^2+5x^3+14x^4+42x^5+132x^6+429x^7+1430x^8+\cdots.
\]
We aim to relate the generating function for simple $321$-avoiding permutations of length at least $4$, which we label $s(x)$, to $c(x)$ and then solve for $s(x)$.

The $321$-avoiding permutations can be divided into four categories:
\begin{itemize}
\item the permutation $1$,
\item the skew decomposable permutations,
\item the sum decomposable permutations, and
\item the inflations of simple permutations of length at least $4$.
\end{itemize}
Obviously the first category of permutations is counted by the generating function $x$.  The skew decomposable $321$-avoiding permutations are all skew sums of two nonempty increasing permutations, and thus are counted by $x^2/(1-x)^2$.  Let $f_\oplus$ denote the generating function for the sum decomposable $321$-avoiding permutations.  Because the $321$-avoiding permutations form a sum closed class ($\pi\oplus\sigma$ avoids $321$ whenever $\pi$ and $\sigma$ both do), we can decompose every sum decomposable $321$-avoiding permutation as the sum of a sum indecomposable permutation and another $321$-avoiding permutation, obtaining
$$
f_\oplus=(c-f_\oplus)c,
$$
which shows that
$$
f_\oplus=\frac{c^2}{1+c}.
$$
Finally, every entry of a simple permutation of length at least $4$ must be involved in an inversion (as otherwise the permutation would be sum decomposable), so to form a $321$-avoiding permutation by inflating such a simple permutation we can only inflate by increasing subsequences.  Thus we see that the contribution of these inflations is $s(x/(1-x))$.  Putting this together shows that
$$
c(x)=x+\frac{x^2}{(1-x)^2}+\frac{c^2}{1+c}+s\left(\frac{x}{1-x}\right),
$$
so we get that
$$
s\left(\frac{x}{1-x}\right)
=
\frac{(1-3x+2x^2-x^3)c-x+x^2-x^3}{(1-x)^2(1+c)},
$$
from which it follows that
$$
s(x)=\frac{1-x-2x^2-2x^3-\sqrt {1-2x-3x^2}}{2+2x}.
$$
The power series expansion of $s(x)$ begins
\[
2 x^4 + 2 x^5 + 7 x^6 + 14 x^7 + 37 x^8 + \cdots,
\]
sequence \OEISlink{A187306} in the \OEISref.

\section{An Iterated System for $321$-Avoiding Simples}

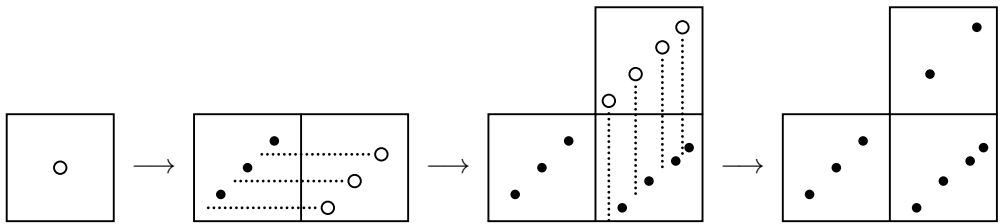
\begin{figure}
\begin{center}
\psset{xunit=0.007in, yunit=0.007in, linewidth=0.01in, runit=0.025in, dotsep=0.015in}
\begin{pspicture}(0,0)(740,160)
	\rput(0,0){
		\psline{c-c}(0,0)(80,0)(80,80)(0,80)(0,0)
		\pscircle(40,40){1.5}
	}
	\rput[c](110,40){$\longrightarrow$}
	\rput(140,0){
		\psline{c-c}(0,0)(160,0)(160,80)(0,80)(0,0)
		\psline{c-c}(80,0)(80,80)
		\psline[linestyle=dotted,linewidth=0.015in](100,10)(10,10)
		\psline[linestyle=dotted,linewidth=0.015in](120,30)(30,30)
		\psline[linestyle=dotted,linewidth=0.015in](140,50)(50,50)
		\pscircle*(20,20){1}
		\pscircle*(40,40){1}
		\pscircle*(60,60){1}
		\pscircle*[linecolor=white](100,10){2}
		\pscircle*[linecolor=white](120,30){2}
		\pscircle*[linecolor=white](140,50){2}
		\pscircle(100,10){1.5}
		\pscircle(120,30){1.5}
		\pscircle(140,50){1.5}
	}
	\rput[c](330,40){$\longrightarrow$}
	\rput(360,0){
		\psline{c-c}(0,0)(160,0)(160,160)(80,160)(80,80)(0,80)(0,0)
		\psline{c-c}(80,0)(80,80)
		\psline{c-c}(80,80)(160,80)
		\psline[linestyle=dotted,linewidth=0.015in](90,90)(90,0)
		\psline[linestyle=dotted,linewidth=0.015in](110,110)(110,20)
		\psline[linestyle=dotted,linewidth=0.015in](130,130)(130,40)
		\psline[linestyle=dotted,linewidth=0.015in](145,145)(145,50)
		\pscircle*(20,20){1}
		\pscircle*(40,40){1}
		\pscircle*(60,60){1}
		\pscircle*(100,10){1}
		\pscircle*(120,30){1}
		\pscircle*(140,45){1}
		\pscircle*(150,55){1}
		\pscircle*[linecolor=white](90,90){2}
		\pscircle*[linecolor=white](110,110){2}
		\pscircle*[linecolor=white](130,130){2}
		\pscircle*[linecolor=white](145,145){2}
		\pscircle(90,90){1.5}
		\pscircle(110,110){1.5}
		\pscircle(130,130){1.5}
		\pscircle(145,145){1.5}
	}
	\rput[c](550,40){$\longrightarrow$}
	\rput(580,0){
		\psline{c-c}(0,0)(160,0)(160,160)(80,160)(80,80)(0,80)(0,0)
		\psline{c-c}(80,0)(80,80)
		\psline{c-c}(80,80)(160,80)
		\pscircle*(20,20){1}
		\pscircle*(40,40){1}
		\pscircle*(60,60){1}
		\pscircle*(100,10){1}
		\pscircle*(120,30){1}
		\pscircle*(140,45){1}
		\pscircle*(150,55){1}
		\pscircle*(110,110){1}
		\pscircle*(145,145){1}
	}
\end{pspicture}
\end{center}
\caption{The development of the simple permutation 2 4 7 1 8 3 5 9 6. In this development, hollow dots might split into groups of entries. In the first step this occurs, as the original single dot splits in three. That imposes two mandatory interpositions (the second and third dots in the second cell). There is one optional addition below it (this addition is, in fact, compulsory in the first step, but optional thereafter). In the next stage of the development only the last dot splits into two, resulting in a single mandatory interposition. One optional addition is also chosen.}
\label{Fig_staircase_simple}
\end{figure}

Every $321$-avoiding permutation $\pi$ has a \emph{staircase decomposition}, as illustrated in the final pane of Figure~\ref{Fig_staircase_simple}.  The square regions of this drawing are called its \emph{cells}. For definiteness we take the elements in the first cell to be the maximum increasing prefix $\tau$ of $\pi$, those in the second cell to be the maximum increasing sequence of values in $\pi \setminus \tau$, and thereafter continue according to the same rules.

We view simple $321$-avoiding permutations as developing one cell at a time, from an initial seed which is a single point. At an intermediate step of this development, the elements in the final cell will represent either single elements of the final permutation, or groups of such elements which are only separated from one another by elements of the next cell. The development of a particular simple permutation is illustrated in Figure \ref{Fig_staircase_simple}.

The discussion above suggests an iterated system which describes the development of $321$-avoiding simple permutations. Specifically, suppose that we have a generating function $s_n(x,y)$ which enumerates all possible developments through $n$ cells, where the elements of the first $n-1$ cells are weighted by $x$ and those of the last cell by $y$. In the next step of the development, each ``hollow dot'' (i.e.~each $y$) can be replaced by a sequence of one or more ``filled dots'' (elements). In the next cell, we \emph{may} put a hollow dot beneath the first of these, and \emph{must} put one between every pair of them. That is, $y$ is replaced by 
\[
x + xy + x^2 y + x^2 y^2 + \cdots = \frac{x(y+1)}{1-xy} .
\]
This gives us the system
\begin{eqnarray*}
s_1(x,y) &=& y, \\
s_{n+1}(x,y) &=& s_n\left(x,  \frac{x(y+1)}{1-xy}\right), \quad \mbox{for $n \geq 1$}.
\end{eqnarray*}

We are interested in the limit as $n \to \infty$ of this system. If we set
\begin{equation}
\label{Eq_FixedPoint}
y = \frac{x(y+1)}{1-xy}\tag{$\dagger$}
\end{equation}
and apply the iteration then we are at a fixed point, and so have produced the desired solution. This yields
\[
s(x) = \frac {1-x-\sqrt {1-2x-3x^2}}{2x},
\]
the power series for the Motzkin numbers,
\[
x + x^2 + 2 x^3 + 4 x^4 + 9 x^5 + 21x^6+51x^7+127x^8+\cdots.
\]
However, the $321$-avoiding simple permutations are not enumerated by the Motzkin numbers. What has gone wrong?

\section{Restricting to Simples}

The development procedure described above fails to generate the simple $321$-avoiding permutations  correctly for two reasons.
\begin{itemize}
\item
The addition of an element below the first dot is optional in this system, but a simple permutation cannot begin with its minimum element, so in the first step this addition should be compulsory.
\item
The addition of an element in the third cell to the left of the absolute minimum element (which is in the second cell) is allowed by this system, but this element (if added) precedes the first descent and thus should be in the first cell.
\end{itemize}
It might seem that these problems would recur in later cells but they do not. In a later horizontal step there is no need to add an element below the smallest element of the preceding cell, since such elements already exist in earlier cells. Similarly, in a vertical step there is no need to forbid the addition of an element to the left of the minimum element, since it is already separated horizontally from the previous top step by at least one element of the previous bottom step. For example, the least element of the fifth cell is separated from the greatest element in the third cell by at least one element of the second cell.

We can correct both these problems by adjusting the initial conditions. One way to do this is to start with a two cell system, and use a new variable $z$ to code the least element. This gives
\[
s_2(x,y,z) = \frac{x z}{1 - xy}.
\]
We then obtain $s_3$ by the substitutions
\[
z \rightarrow \frac{x}{1-xy} \quad\mbox{and}\quad y \rightarrow \frac{x(y+1)}{1-xy}.
\]
This eliminates $z$ and thereafter we use the original iteration. It follows that we should obtain the generating function for simple $321$-avoiding permutations by substituting (\ref{Eq_FixedPoint}) into $s_3$. Doing so yields
\[
s(x) = \frac{1-x-\sqrt {1-2x-3x^2}}{2+2x}
\]
whose power series expansion begins
\[
x^2 + 2 x^4 + 2 x^5 + 7 x^6 + 14 x^7 + 37 x^8 + \cdots,
\]
which (once the term corresponding to the permutation $21$ is removed) agrees with the result obtained in Section~\ref{sec-implicit-321}.

\section{Skew-Merged Permutations}

\begin{figure}
\centerline{
\psset{xunit=0.01in, yunit=0.01in} \psset{linewidth=0.005in}
\begin{pspicture}(0,-10)(170,170)
\multips{0}(0,0)(60,60){3}{\pspolygon[linestyle=solid](0,0)(50,0)(50,50)(0,50)}
\multips{0}(0,120)(120,-120){2}{\pspolygon[linestyle=solid](0,0)(50,0)(50,50)(0,50)}
\psline(5,5)(45,45)
\psline(125,45)(165,5)
\psline(125,125)(165,165)
\psline(5,165)(45,125)
\rput(10,40){I}
\rput(10,130){IV}
\rput(160,130){III}
\rput(160,40){II}
\end{pspicture}
}
\caption{The basic structure of a skew-merged permutation as described in \cite{atkinson:permutations-wh:}.  Area I consists of all the elements that play the role of $1$ in some $132$, area II those that play the role of $1$ in some $231$, area III those that play the role of $3$ in some $213$ and area IV those that play the role of $3$ in some $312$. In general, any or all of these areas might be empty, and each is monotone of the type indicated by the line segments. The central region consists of all the remaining elements and could be either increasing or decreasing.}
\label{Fig_skew_merged}
\end{figure}
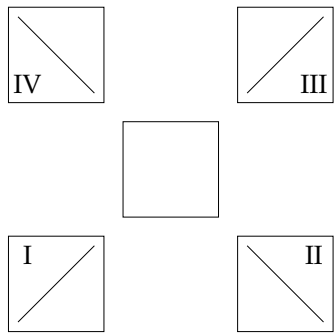

Recall that a skew-merged permutation is one that can be written as the union of an increasing and a decreasing sequence. The basic structure of such a permutation is shown in Figure \ref{Fig_skew_merged}. We intend to apply an analogous technique to that of the previous sections to enumerate the simple skew-merged permutations. In order to do so, we must refine that basic structure somewhat. The first thing to note is that the central area is an interval, and hence in a simple permutation can contain at most one point. Every simple permutation of length at least four contains either $2413$ or $3142$ as a subpermutation, so if $\pi$ is a simple skew-merged permutation, then all four areas are occupied by at least one point. Further, it is easy to check that if $\pi$ has a central element, $c$, and is simple, then $\pi-c$ (the permutation obtained by $\pi$ by removing $c$ and relabeling the remaining entries) is also simple. We call the four elements (one from each area) closest to the centre of $\pi$ its \emph{inner elements}. 

To continue we assume that $\pi$ has no central element and further divide into two (symmetric) cases depending on the pattern of the four inner elements. There are only two possible patterns of this type, and the following discussion assumes that we are considering a simple skew-merged permutation $\pi$ whose inner elements have the pattern $3142$ (the other possibility is its inverse, $2413$). Denote the inner elements specifically as $cadb$ from left to right.

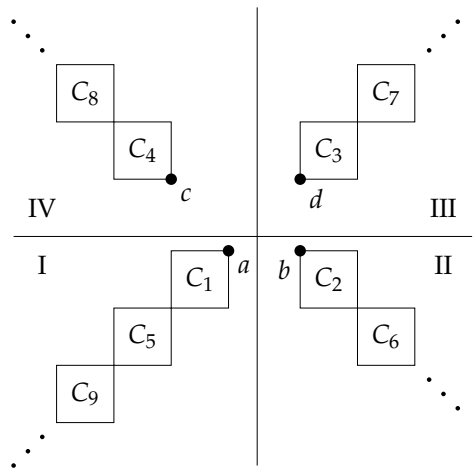
\begin{figure}
\centerline{
\psset{xunit=0.015in, yunit=0.015in, runit=0.015in} \psset{linewidth=0.005in}
\begin{pspicture}(0,-10)(160,160)
\psline(85,0)(85,160)
\psline(0,80)(160,80)
\multips{0}(15,15)(20,20){3}{\pspolygon[linestyle=solid](0,0)(20,0)(20,20)(0,20)}
\multips{0}(15,120)(20,-20){2}{\pspolygon[linestyle=solid](0,0)(20,0)(20,20)(0,20)}
\multips{0}(120,120)(-20,-20){2}{\pspolygon[linestyle=solid](0,0)(20,0)(20,20)(0,20)}
\multips{0}(120,35)(-20,20){2}{\pspolygon[linestyle=solid](0,0)(20,0)(20,20)(0,20)}
\pscircle*(75,75){2}
\uput[-45](75,75){$a$}
\pscircle*(100,75){2}
\uput[-135](100,75){$b$}
\pscircle*(100,100){2}
\uput[-45](100,100){$d$}
\pscircle*(55,100){2}
\uput[-45](55,100){$c$}
\rput(25,25){$C_9$}
\rput(25,130){$C_8$}
\rput(130,130){$C_7$}
\rput(130,45){$C_6$}
\rput(45,45){$C_5$}
\rput(45,110){$C_4$}
\rput(110,110){$C_3$}
\rput(110,65){$C_2$}
\rput(65,65){$C_1$}
\rput(10,70){I}
\rput(10,90){IV}
\rput(150,90){III}
\rput(150,70){II}
\multips{0}(0,0)(5,5){3}{\pscircle*(0,0){0.8}}
\multips{0}(145,145)(5,5){3}{\pscircle*(0,0){0.8}}
\multips{0}(0,155)(5,-5){3}{\pscircle*(0,0){0.8}}
\multips{0}(145,30)(5,-5){3}{\pscircle*(0,0){0.8}}
\end{pspicture}
}
\caption{The decomposition of a simple skew-merged permutation whose inner non-central elements $cadb$ have pattern $3142$ into monotone cells, each interacting with only its immediate predecessor and successor.}
\label{Fig_skew_simple}
\end{figure}

Now we describe a partition of $\pi$ into a sequence of monotone cells (akin to the staircase decomposition of $321$-avoiding permutations) where each cell in the sequence interacts only with its immediate predecessor and its immediate successor. These cells are arranged in a counterclockwise spiral pattern, moving outwards. Their description is most easily understood by referring to Figure \ref{Fig_skew_simple}. The first four of these cells  are defined specifically as follows.
\begin{itemize}
\item
$C_1$ consists of all the elements in area I lying to the right of $c$. In particular $a \in C_1$.
\item
$C_2$ consists of all the elements in area II which lie above some element of $C_1$. In particular, $b \in C_2$.
\item
$C_3$ consists of all the elements in area III which lie to the left of some element of $C_2$. In particular, $d \in C_3$.
\item
$C_4$ consists of all the elements in area IV which lie below some element of $C_3$. In particular, $c \in C_4$. Note also that $C_4$ lies entirely to the left of $C_1$ by the initial choice of $C_1$.
\end{itemize}

Now for $k \geq 1$ the cells are defined inductively as follows.
\begin{itemize}
\item
$C_{4k+1}$ consists of all the elements in area I which lie to the right of some element of $C_{4k}$ but do not belong to some previous cell.
\item
$C_{4k+2}$ consists of all the elements in area II which lie above some element $C_{4k+1}$ but do not belong to some previous cell.
\item
$C_{4k+3}$ consists of all the elements in area III which lie lie to the left of some element of $C_{4k+2}$ but do not belong to some previous cell.
\item
$C_{4k+4}$ consists of all the elements in area IV which lie below some element of $C_{4k+3}$ but do not belong to some previous cell.
\end{itemize}

If $C_n$ is empty for some $n$ then the elements of $\pi$ that belong to the union of the $C_k$ for $k < n$ form an interval. Since $\pi$ is assumed to be simple, this must be the entire permutation.

Now we reverse the perspective above, which deconstructs a simple skew-merged permutation $\pi$ into a sequence of cells and view it as a constructive recipe for building these permutations. As before we think of the cells $C_1$ through $C_{k-1}$ as having been constructed for some $k$, along with some elements of $C_k$ represented by hollow dots. In extending the construction, these hollow dots can split into sequences which then impose mandatory insertions in the next cell, while the spaces between hollow dots (or before the first one in a cell) represent optional insertions. Note that the precise meaning of ``before'' depends on the area in which we are constructing the next cell, for example, in area II such an element would be the leftmost and greatest element of $C_{k+1}$. We code the generating function for the result of applying this iteration through cell $n$ by $u_n(x,y)$ where $y$ labels the hollow dots.
Because of the necessity of including the four extreme points in the first four cells, the first three optional insertions at the beginning of a cell are in fact compulsory. The fourth one is actually forbidden (it would violate the ``inner'' property of $c$). We account for these variations to the general rule by adjusting the substitutions corresponding to these steps:
$$
\begin{array}{lcll}
u_1(x,y) &=& y, \\
u_2(x,y) &=& u_1\displaystyle\left(x, \frac{x(1+y)}{1-xy} \right) \cdot\frac{y}{y+1}, \\[12pt]
u_3(x,y) &=& u_2\displaystyle\left(x, \frac{x(1+y)}{1-xy} \right) \cdot\frac{y}{y+1}, \\[12pt]
u_4(x,y) &=& u_3\displaystyle\left(x, \frac{x(1+y)}{1-xy} \right) \cdot\frac{y}{y+1}, \\[12pt]
u_5(x,y) &=& u_4\displaystyle\left(x, \frac{x(1+y)}{1-xy} \right) \cdot\frac{1}{y+1}, \\[12pt]
u_{n+1}(x,y) &=& u_n\displaystyle\left(x, \frac{x(1+y)}{1-xy}\right),&\mbox{for $n \geq 5$}.
\end{array}
$$
By computing $u_5$ and then substituting from (\ref{Eq_FixedPoint}) as before, we obtain the generating function for the simple skew-merged permutations of this type:
\[
u = \frac{1 - 2x - x^2  + (x-1) \sqrt{1 - 2x - 3x^2}}{2(x+1)^2}.
\]
The generating function for all simple permutations of length at least $4$ in this class is then $s=2(x+1)u$.  Here the $2$ accounts for the two types of inner pattern and the $x+1$ for the possible addition of a central element.  The power series expansion of $s(x)$ begins
\[
2 x^4 + 2 x^5 + 8 x^6 + 16 x^7 + 44 x^8 + \cdots,
\]
sequence \OEISlink{A220589} in the \OEISref.

Finally, we can recover the generating function, $f$, for the class of skew-merged permutations from its simple elements by the usual techniques.  Skew-merged permutations are either sum or skew decomposable or have a simple quotient of length at least $4$. Every sum decomposable permutation in this class is of the form $1 \oplus \pi$ or $\pi \oplus 1$ for a skew-merged permutation $\pi$, and the skew decomposable permutation can be analogously described. Allowing for the over-counting of elements of the form $1 \oplus \pi \oplus 1$ or $1 \ominus \pi \ominus 1$ we see that the sum or skew decomposable elements of $\S$ are enumerated by
\[
4xf - 2x^2(f+1).
\]
In an element with a simple quotient of length at least $4$, if there is no central element in the quotient then the quotient can only be inflated by monotone sequences at each point. If a central element is present then it can be inflated by an arbitrary skew-merged permutation. Combining this with our previous observations gives:
\[
f = x + 4 x f - 2 x^2 (f+1) + 2 u \left( \frac{x}{1-x} \right) (f + 1).
\]
And, although it is hardly evident at a glance, this gives the known result
\[
f
=
\frac{1-3x}{(1-2x)\sqrt{1-4x}}
=
1+x+2x^2+6x^3+22 x^4 + 86 x^5 + 340 x^6 + 1340 x^7 + 5254 x^8 + \cdots.
\]
\OEIS{A029759}  To confirm this, we suggest that the reader do as we did and leave the details to a computer algebra system, but the main reason that things simplify so nicely is that the substitution of $x/(1-x)$ for $x$ in $\sqrt{1-2x-3x^2}$ yields $\sqrt{1-4x}/(1-x)$.
\section{Extensions}


It has not escaped our notice that the specific techniques we have detailed suggest possible extensions.  In the language of ``geometric grid classes'' studied in \cite{albert:geometric-grid-:,albert:inflations-of-g:}, the $321$-avoiding permutations can be described as
$$
\Geom
\begin{footnotesize}
\left(
\begin{array}{ccccccc}
&&&&&\reflectbox{$\ddots$}&\reflectbox{$\ddots$}\\
&&&&1&1\\
&&&1&1\\
&&1&1\\
&1&1\\
1&1\\
\end{array}
\right),
\end{footnotesize}
$$
and we have shown that the skew-merged permutations can be expressed as
$$
\Geom
\begin{footnotesize}
\left(
\begin{array}{rrrrrrr}
\ddots&&&&&&\reflectbox{$\ddots$}\\
&-1&&&&1\\
&&-1&&1\\
&&&1&-1\\
&&1&&&-1\\
&\reflectbox{$\ddots$}&&&&&\ddots
\end{array}
\right)
\end{footnotesize}
\bigcup
\Geom
\begin{footnotesize}
\left(
\begin{array}{rrrrrrr}
&\ddots&&&&&\reflectbox{$\ddots$}\\
&&-1&&&1\\
&&&-1&1\\
&&&1\\
&&1&&-1\\
&1&&&&-1\\
\reflectbox{$\ddots$}&&&&&&\ddots
\end{array}
\right).
\end{footnotesize}
$$
It should be possible to use our techniques to enumerate the simple permutations for all infinite geometric grid classes which consist of a single ``path'', so long as this path satisfies some suitable regularity condition.  Of course, the analysis of which hollow dots are mandatory and which are optional will depend on the specific grid class studied.

\bibliographystyle{acm}
\bibliography{../refs}

\end{document}